\documentclass{amsart}

\usepackage{amsmath}
\usepackage{amssymb}
\usepackage{amscd}
\usepackage{amsthm}
\usepackage{tikz}
\usepackage{tikz-cd}
\usepackage{graphicx}
\usepackage[pagebackref]{hyperref}

\newcommand{\lra}{\longrightarrow}

\newcommand{\rsp}{\raisebox{0em}[2.7ex][1.3ex]{\rule{0em}{2ex} }}

\newcommand{\Z}{{\mathbb Z}}

\newcommand{\eps}{\varepsilon}

\newcommand{\cM}{{\operatorname{\mathcal M}}}

\title{Gauss and $p$-adic numbers}
\author{F. Lemmermeyer}

\begin{document}

\maketitle

\begin{abstract}
  Gauss's notebooks contain endless calculations that show how he
  applied known algorithms or developed them on the fly. In this
  article we will take a look at his toolkit and explain most of
  his ideas concerning calculations with ``infinite congruences'',
  which is Gauss's name for what Hensel later called p-adic numbers.  
\end{abstract}

The notion of $p$-adic numbers is due to Hensel, who introduced them
in 1899 (see \cite{Hensel}). Hasse made them an indispensable tool in
algebraic number theory by his discovery of Local-Global principles,
at first in connection with quadratic forms, later in the theory of
algebras and norm residues.

It is also well known that Gauss knew some form of ``Hensel's Lemma'',
which he used in the planned Section VIII of his Disquisitiones
(see Frei's article in \cite{GSS}).
It turns out that Gauss had played around with ``infinite congruences''
at the same time; actually, these infinite congruences modulo $p^\infty$
are $p$-adic numbers! In his notebook \cite{GS}  titled ``1800 Juli'',
Gauss recorded several calculations with these infinite congruences
modulo $241$, modulo $11$ and modulo $10$. In this article we will explain
most of Gauss's calculations and leave the problem of figuring out
the rest as a challenge to the readers.

\section{Square roots and cube roots}

In his notebooks, Gauss used different methods for computing roots.
In this section we will discuss a few of them.

\subsection*{Square roots using the binomial theorem}

The binomial expansion of $\sqrt{1 \pm x}$ (this is just the Taylor expansion
around $x = 0$) is given by
\begin{align}
  \label{BE-} \sqrt{1-x} & = 1 - \frac12 x - \frac18x^2
        - \frac1{16}x^3 - \frac{5}{128}x^4
        - \frac{7}{256} x^5 - \cdots, \\
  \label{BE+} \sqrt{1+x} & =   1 + \frac12 x - \frac18x^2
        + \frac1{16}x^3 - \frac{5}{128}x^4
        + \frac{7}{256} x^5 - \cdots.
\end{align}
        
For computing $\sqrt{5}$ we could use the fact that $5 = 4 + 1$, which implies
$$ \sqrt{5} = \sqrt{4+1} = 2 \sqrt{1 + \frac14}  =
      2\Big(1 + \frac12 \cdot \frac14 - \frac18 \cdot \frac1{16} +  \ldots \Big)
      \approx \frac{143}{64}. $$
We get a series that converges faster if we  use a number of the form $5n^2$ 
close to a square; if we choose $5 \cdot 4^2 = 80 = 9^2-1$, then we find
$$ 4\sqrt{5}  = \sqrt{80} = \sqrt{9^2-1} = 9 \sqrt{1 - \frac1{81}}. $$
Applying the binomial expansion (\ref{BE-}) we get
$$ 4 \sqrt{5}
  = 9 \Big(1 - \frac1{2 \cdot 81} - \frac1{8 \cdot 81^2} - \ldots \Big), $$
which yields
$$ \sqrt{5} \approx \frac{52163}{23328} \approx 2.2360682. $$

\subsubsection*{Computing $\sqrt{79}$}

In his investigations (see \cite[p.~122]{GB}) on algebraic numbers of the
form $a + b\sqrt{m} + c \sqrt{n} + d \sqrt{mn}$, Gauss needs an approximation
of $\sqrt{79}$; in his calculation, the only intermediate results he writes
down are the digits 7 1 1 0 5. The result he obtains is
$$ \sqrt{79} \approx 8.8882. $$
I do not know which algorithm Gauss used.

\begin{figure}[ht!]
  \begin{center}
    \includegraphics[width=6cm]{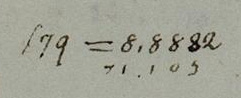}
  \end{center}
  \caption{Computation of $\sqrt{79}$}
\end{figure}

\subsubsection*{$\sqrt{79}$ via binomial expansion.}
We can compute an approximation of $\sqrt{79}$ using
$$\sqrt{79} = \sqrt{9^2-2} = 9 \sqrt{1 - \frac2{81}}; $$
This shows
\begin{align*}
  \sqrt{79}  & \approx
  9\Big(1 - \frac12 \cdot \frac2{81} - \frac18 \cdot \frac4{81^2} \Big) \\
  & = 9 - \frac19 - \frac1{2 \cdot 3^6} \approx 8.8882
\end{align*}
Since $9 - \frac19 = 8.88888\ldots$, Gauss only needed to subtract
$\frac1{2 \cdot 729}$.

We can also compute $\sqrt{79}$ quickly using continued fractions. We have
$$ \sqrt{79} = 8 +
  \cfrac{1}{1+\cfrac{1}{7+\cfrac{1}{1+\cfrac{1}{16+\cfrac{1}{1+\cdots}}}}}. $$
The corresponding fractions approximating $\sqrt{79}$ are
$$ \frac{71}8, \quad \frac{80}9, \quad \frac{1351}{152}, \quad
   \text{and} \quad \frac{1431}{161}. $$
The first fraction would explain Gauss's $71$, but I cannot
see how the remaining digits are connected to the continued fraction.

\subsection*{Square Root Algorithm}
Next we will explain the lost art of extracting square roots digit by
digit.  For computing the square root of $N = 717\,409$, we write this
number as $71|74|09$; then $8^2 < 71 < 9^2$ tells us that
$800 < \sqrt{N} < 900$. Thus the leftmost digit of the square root is
$a = 8$, which we write on the left and on the top of the square root:

\begin{center}
  \begin{tikzpicture}
    \draw (0,0) node  {$\sqrt{71|74|09}$};
    \draw (-1.3,0) node {$8$};
    \draw (0, 0.5) node {$8$};
  \end{tikzpicture}
\end{center}

We multiply the $8$ on top with the $8$ on the left, subtract
the product $64$ from $71$, and append the next two digits $74$
to the difference $7$. Double the $8$ on top and append a $0$; then
$160$ goes $4$ times into $774$; we append the $4$ to the $8$ on top
and to $2 \cdot 8 = 16$ on the left and subtract $4 \cdot 164 = 656$
from $774$. We append the final digits $09$ to the difference; double
the $84$ on top; $168$ goes $7$ times into $1180$.
Append a $7$ to the $84$ on top and to $168$ on the left, and subtract
$7 \cdot 1687 = 11809$. The result is $0$, and $\sqrt{717409} = 847$.

\begin{center}
  \begin{tikzpicture}
    \draw (0,0) node  {$\sqrt{71|74|09}$};
    \draw (-1.3,0) node[anchor=east] {$8$};
    \draw (0, 0.5) node {$8 \ 4 \ 7$};
    \draw (-0.22,-0.5) node {$64|$};
    \draw (0.09,-0.9) node   {$ 7|74$};  
    \draw (-1.3,-0.9) node[anchor=east] {$164$};
    \draw (0.09,-1.3) node {$ 6|56$};
    \draw (0.3,-1.7)   node {$ 1|18|09$}; 
    \draw (-1.3,-1.7) node[anchor=east] {$1687$};  
    \draw (0.3,-2.1)   node {$ 1|18|09$};   
  \end{tikzpicture}
\end{center}

Let us now explain why this algorithm works. The first step is clear:
We choose our first digit $a$ in such a way that $a^2 < 71 < (a+1)^2$,
i.e., we choose $a = 8$. Next we need to find the maximal digit $b$
with $(80+b)^2 < 7174$; this is equivalent to $(160 + b)b < 7174 - 6400 = 774$.
Since $160$ goes four times into $774$, we try $b = 4$ (occasionally
we overshoot and then have to decrease our choice by $1$).

Now we try to find a digit $c$ with $(840+c)^2< 717\,409$. This
is equivalent to $(1680 + c)c < 717\,409 - 840^2 = 11\,809$. We have
obtained the remainder by first subtracting $800^2$ and then
$16400 \cdot 4$. This is in line with
$$ N - (a+b)^2 = N^2 - a^2 - (2a+b)b, $$
which in our case, where $N = 717\,409$, $a = 800$ and $b = 40$, yields
$$ N - 840^2 = N - 800^2 - 16400 \cdot 4. $$ 
Since $1680$ goes seven times into $11809$ we choose $c = 7$, and now
$$ N - 847^2 = N - 840^2 - 2 \cdot 840 \cdot 7 - 7^2
             =  11809 - 1687 \cdot 7 = 0. $$
With a little practice, this algorithm can be performed easily on
numbers of moderate size.  The calculation of $\sqrt{79}$ using this
method would look as follows.

\begin{center}
\begin{tikzpicture}
   \draw (0,0) node[anchor=east] {$\sqrt{79}\phantom{00|00|00}$};
    \draw (-2.5,0) node[anchor=east] {$8$};
    \draw (-1, 0.5) node {$8 \ 8 \  8  \ 8  \ 2$};
    \draw (-0,-0.5) node[anchor=east] {$64|\phantom{00|00|00}$};
    \draw ( 0,-.9) node[anchor=east]    {$15|00\phantom{|00|00}$};  
    \draw (-2.5,-.9) node[anchor=east] {$168$};
    \draw ( 0,-1.3) node[anchor=east] {$13|44\phantom{|00|00}$};
    \draw ( 0,-1.7)  node[anchor=east] {$ 1|56|00\phantom{|00}$}; 
    \draw (-2.5,-1.7) node[anchor=east] {$1768$};  
    \draw (0,-2.1)  node[anchor=east] {$1|41|44\phantom{|00}$};  
    \draw (0,-2.5)  node[anchor=east] {$14|56|00$};    
    \draw (-2.5,-2.5) node[anchor=east] {$17768$};   
    \draw (0,-2.9)  node[anchor=east] {$14|21|44$}; 
    \draw (0,-3.3)  node[anchor=east] {$34|56$};
    \draw(-2.5,-3.3) node[anchor=east] {$177762$};
\end{tikzpicture}
\end{center}
This approximation is slightly too large.

\subsection*{Cube roots}

Gauss had a lot of algorithms he could choose from when doing calculations.
In his notes on units of the form $a + b\sqrt{2} + c \sqrt{4}$
(see \cite[p.~161]{GB}) he computes an approximation of $\sqrt[3]{2}$
using the inequality of arithmetic and geometric mean.

The corresponding algorithm for square roots is straightforward. Starting
with $a = \frac32$ choose $b = \frac2{a} = \frac43$; then
$$ \sqrt{2} = \sqrt{ab} \le \frac{a+b}2 = \frac{17}{12}. $$
Repeating this step with $a = \frac{17}{12}$ and $b = \frac2a = \frac{24}{17}$
yields
$$ \sqrt{2} \approx \frac{577}{408} \approx 1.414215. $$

For the cube root of $2$, Gauss's calculations are a little bit roundabout.
Apparently Gauss uses the elementary inequalities
$$ \frac{m+1}{n+1} < \frac{m}{n} < \frac{m-1}{n-1} $$
for fractions with $m > n > 1$.

Gauss starts with $a = b = 1$ and $c = \frac2{ab} = 2$ and gets
$$ \sqrt[3]{2} = \sqrt[3]{abc}  \le \frac{a+b+c}3 = \frac43. $$
He writes down an ``error term'' $+0.073 \approx \frac43 - \sqrt[3]{2}$ which
I think he computed at the end.

\begin{figure}[ht!]
  \begin{center}
    \includegraphics[width=9cm]{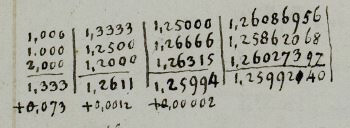}
  \end{center}
  \caption{Computation of $\sqrt[3]{2}$ using
    arithmetic-geometric mean}\label{Gsqrt3}
\end{figure}

In the next step, he takes $a = \frac43 > \sqrt[3]{2}$; he chooses a
fraction that is slightly smaller, namely $b = \frac{4+1}{3+1} = \frac54$,
and gets $c = \frac2{ab} = \frac65$. Applying the inequality between
geometric and arithmetic mean again he finds
$$ \sqrt[3]{2} = \sqrt[3]{abc}  \le \frac{a+b+c}3  =
   \frac{227}{180} \approx 1.2611. $$

Next Gauss picks the fraction among $a, b, c$ closest to the new
approximation, namely $\frac54$, which is slightly too small, and
looks for a fraction $> \frac54$ with denominator divisible by $5$. Since
$\frac54 = \frac{20}{16}$, he picks
$b = \frac{20-1}{16-1} = \frac{19}{15} \approx 1.26666$ and gets
$c = \frac2{ab} = \frac{24}{19} \approx 1.26315$, and now
$$  \sqrt[3]{2}  \approx   \frac{4309}{3420} \approx 1.25994. $$

In his final step, Gauss may have used the continued fraction algorithm
for his approximation $\sqrt[3]{2} \approx 1.25994$, namely
$$ 1.25994 \approx 
   1 + \cfrac{1}{3 + \cfrac{1}{1+\cfrac{1}{5}}} = \frac{29}{23}. $$
Looking for a fraction close to $\frac54$ with denominator divisible
by $29$ he finds $b = \frac{75-2}{60-2} = \frac{73}{58}$.
Thus in the last row he computes the arithmetic mean of
$$ a = \frac{29}{23} = \frac{5 \cdot 6 -1}{4 \cdot 6 - 1}, \quad  
   b = \frac{73}{58} \quad \text{and} \quad c = \frac{92}{73}, $$ 
giving him the approximation
$$ \sqrt[3]{2} \approx \frac{368\,081}{292\,146} \approx 1.259\,921\,409. $$

This is the whole calculation as recorded in his notebook
(see Fig.~\ref{Gsqrt3}):
$$ \begin{array}{r|r|r|r}
  \rsp 1 & \frac43 & \frac{ 5}{ 4} & \frac{29}{23} \\
  \rsp 1 & \frac54 & \frac{19}{15} & \frac{73}{58}  \\
  \rsp 2 & \frac65 & \frac{24}{19} & \frac{92}{73}
 \end{array}  \qquad 
 \begin{array}{r|r|r|r}
   \rsp 1.000 & 1.3333 & 1.25000 &  1.26086956 \\
   \rsp 1.000 & 1.2500 & 1.26666 &  1.25862068 \\
   \rsp 2.000 & 1.2000 & 1.25994 &  1.26027397
 \end{array} $$

The continued fraction expansion of $\sqrt[3]{2}$ begins with
$$ \sqrt[3]{2} = 1 + \cfrac{1}{3 + \cfrac{1}{1+\cfrac{1}{5+
                \cfrac{1}{1 + \cfrac{1}{1 + \cfrac{1}{4 + \cdots}}}}}}. $$
The first few rational approximations then are
$$ a_1 = 1, \quad a_2 = \frac43, \quad a_3 = \frac54, \quad
   a_4 = \frac{29}{23}, \quad a_5 = \frac{34}{27}, \quad 
   a_6 = \frac{63}{50}, \quad a_7 = \frac{286}{227}. $$
Quite possibly Gauss chooses $b$ using an idea I cannot see.

\section{$p$-adic roots of polynomials}

The editors of Gauss's Werke published several short notes from
Gauss's notebooks; a lot of material was left unpublished and has now
been made accessible online. The main topic of this article concerns
a few pages in Gauss's notebook from July 1800, where Gauss is working
with congruences modulo infinite powers of integers: these are, as we will
see, $g$-adic numbers.

In \cite[p.~14]{GS}, Gauss presents the following calculation
(Fig.~\ref{F241}; observe the impeccable handwriting):
\medskip

\begin{figure}[ht!]
  \includegraphics[width=9cm]{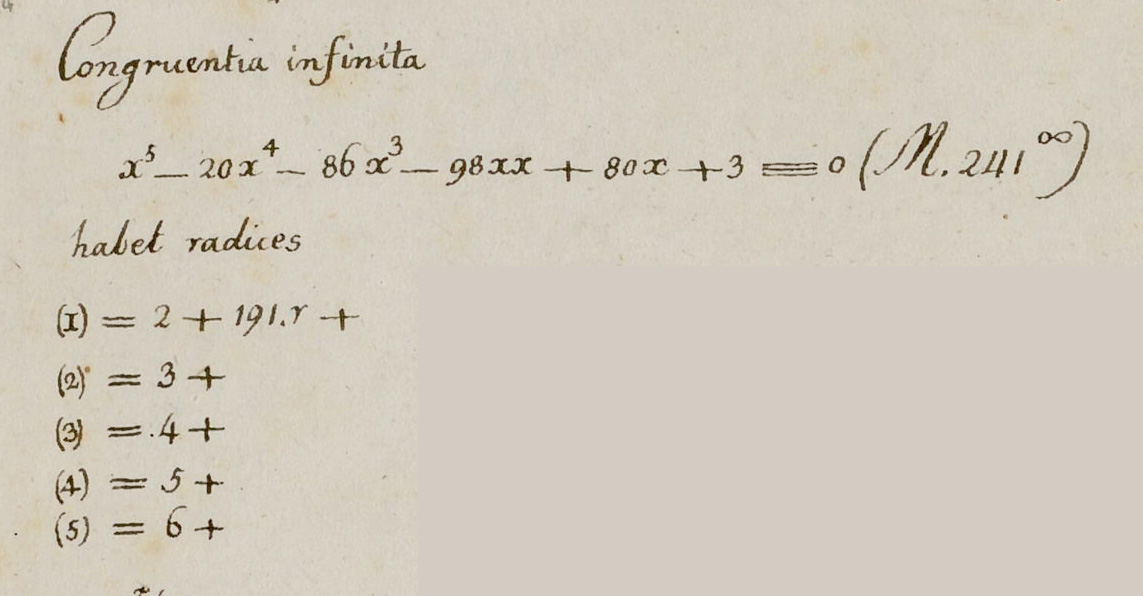}
  \caption{Congruences modulo $241^\infty$}\label{F241}
\end{figure}

\begin{quote}
  {\large {\em Congruentia infinita}}
  $$ x^5 - 20x^4 - 86x^3 - 98x^2 + 80x + 3 \equiv 0 \quad(\cM.241^\infty) $$
  {\em habet radices}
  \begin{align*}
    (1) & = 2 + 191  \cdot r + \\
    (2) & = 3 + \\
    (3) & = 4 + \\
    (4) & = 5 + \\
    (5) & = 6 +
  \end{align*}
\end{quote}

Gauss apparently computed the polynomial
$$ (x-2)(x-3)(x-4)(x-5)(x-6) = x^5 - 20x^4 + 155x^3 - 580x^2 + 1044x - 720 $$
and reduced the coefficients to their smallest values modulo $241$; this
guaranteed that his polynomial
$$ f(x) = x^5 - 20x^4 - 86x^3 - 98x^2 + 80x + 3 $$
has roots $x_1 \equiv 2$, $x_2 \equiv 3$, $x_3 \equiv 4$, $x_4 \equiv 5$
and $x_5 \equiv 6$ modulo $241$. Then he starts computing $p$-adic
approximations of these roots modulo $241^2$; for the first root he finds
$$ x_1 \equiv 2 + 191 \cdot r \bmod 241^2, $$
where apparently $r = 241$. The lifts modulo $241^3$ of these roots are
$$ \begin{array}{cccrcrcl}
  x_1 & = 2 & + & 191 \cdot 241 & + & 160 \cdot 241^2 & + \ldots, \\
  x_2 & = 3 & + & 238 \cdot 241 & + &  16 \cdot 241^2 & + \ldots, \\
  x_3 & = 4 & + & 192 \cdot 241 & + & 221 \cdot 241^2 & + \ldots, \\
  x_4 & = 5 & + &  65 \cdot 241 & + &  17 \cdot 241^2 & + \ldots, \\
  x_5 & = 6 & + &  37 \cdot 241 & + &  65 \cdot 241^2 & + \ldots.
\end{array} $$
Gauss was obviously aware of the fact that all the roots can be lifted
to roots modulo arbitrarily large powers of $p = 241$; in the limit, the
roots modulo $241^\infty$ are the $p$-adic roots of $f$. Gauss did not
compute these approximations; instead he turned to a more interesting
problem.

\section{$p$-adic approximations of quadratic Gauss sums}

At the bottom of the page with the $241$-adic approximations of the
roots of a polynomial, Gauss writes
\begin{quote}
  $$ \sqrt{5} \pmod{11^\infty} = 9.0.4.10.4.4$$  
\end{quote}
This is the $11$-adic expansion of one of the two square roots
modulo $5$, to be read from right to left:
$$ \sqrt{5} = 4 + 4 \cdot 11 + 10 \cdot 11^2 + 4 \cdot 11^3 + 0 \cdot 11^4
                + 9 \cdot 11^5 + \ldots $$
Using {\tt pari}, this approximation is easily computed as
$$ {\tt sqrt(5 + O(11^8))} =
     4 + 4 \cdot 11 + 10 \cdot 11^2 + 4\cdot 11^3 + 9 \cdot 11^5
       + 5 \cdot 11^6 +  8 \cdot 11^7 + O(11^8). $$

       \medskip

\begin{figure}[ht!]
  \includegraphics[width=8cm]{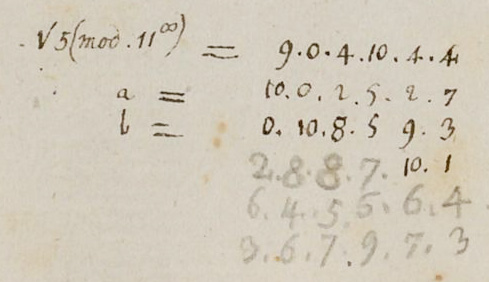}  
  \caption{Square root of $5$ in the $11$-adic numbers}\label{Gab}
\end{figure}

Now let $p$ be an odd prime and $a$ an odd natural number not
divisible by $p$. Solve the congruence $an^2 \equiv 1 \bmod 2p$ and
write $an^2 = 1 + 4pk$ for some integer $k$. Using the binomial expansion
of $\sqrt{1+x}$ we find
\begin{align*}
  \sqrt{an^2}
  & = \sqrt{1 + 4kp} \\
  & = 1 + 2kp - 2k^2p^2 + 4k^3p^3 - 10k^4p^4 + 28k^5p^5 - 84k^6p^6 + \ldots \\
  & = (1 + 2kp + 4k^3p^3 + 28k^5p^5 + \ldots)
      - (2k^2p^2 +  10k^4p^4 + 84k^6p^6 + \ldots).
\end{align*}
In the special case where $p = 11$ and $a = 5$ we have
$a \cdot 3^2 = 45 = 1 + 4 \cdot 11$, hence $k = 1$; since $28 \equiv 6 \bmod 11$
we find
$$ \begin{array}{ccccrrrrrr}
            &   &   &     & 11^5 & 11^4 & 11^3  & 11^2 & 11 & 1 \\  
  \sqrt{45} & = &   & \ldots &   6  &   0 &   4  &  0    & 2   & 1 \\
            &   & - & \ldots &   0  &  10 &   0  &  2    & 0   & 0 \\ \hline
            &   &   & \ldots &   5  &   1 &   3  &  9    & 2   & 1 \\
            &   &   & \ldots &   9  &   0 &   4  & 10    & 4   & 4
\end{array} $$
Here, the second line means that $\sqrt{45}$ is congruent to
$6 \cdot 11^5 + \ldots + 1 \bmod 11^6$. In the last line, we have divided the
result by $3$ in order to obtain $\sqrt{5} = \frac{\sqrt{45}}3$. These are
exactly the calculations performed by Gauss:

\begin{figure}[ht!]
  \includegraphics[width=8cm]{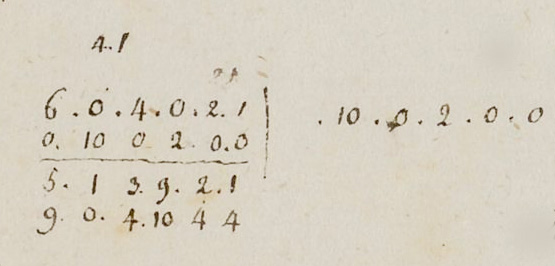}
  \caption{Computation of an $11$-adic approximation of $\sqrt{5}$}\label{Ab5}
\end{figure}

The meaning of the numbers $a$ and $b$ in Fig.~\ref{Gab} is explained
elsewhere on this page:

\begin{figure}[ht!]
  \includegraphics[width=8cm]{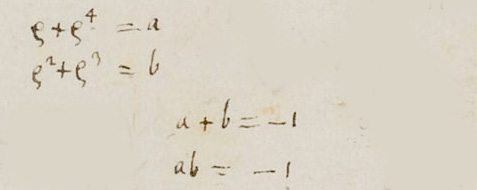}
  \caption{Computation of $11$-adic approximations of quadratic periods}
\end{figure}

Let $\rho$ be a primitive fifth root of unity. Then the quadratic
periods are $a = \rho + \rho^4 = \frac{-1+\sqrt{5}}2$ and 
$b = \rho^2 + \rho^3 = \frac{-1-\sqrt{5}}2$; observe that
$a + b = -1$ and $ab = -1$; the difference
$a - b = \rho - \rho^2 - \rho^3 + \rho^4 = \sqrt{5}$ is the
quadratic Gauss sum\footnote{On \cite[p.~34]{GS} Gauss presents a
  ``theorema novissimum pulcherrimum'', a new beautiful theorem,
  namely the determination of the sign of quadratic Gauss sums.}
modulo $5$.

We find $2a = -1 + \sqrt{5} = 9.0.4.10.4.3$ and therefore
$a =    4.5.7.10.7.7$, as well as $b = 6.5.3. 0.3.3$.
The values given by Gauss apparently contain a mistake:

$$ \begin{array}{ccrrrrrr}
  a & = & 10 &  0 & 2 & 5 & 2 & 7 \\
  b & = &  0 & 10 & 8 & 5 & 9 & 3 \\
    &   &  2 & 8 & 8 & 7 & 10 & 1 \\
    &   &  6 & 4 & 5 & 6 &  6 & 4 \\
    &   &  3 & 6 & 7 & 9 &  7 & 3
\end{array} $$

When dividing $9.0.4.10.4.3$ by $2$, borrowing $1$ from left gives
$3 +11 = 14$, and dividing by $2$ gives the final digit $7$. Gauss
forgot that he has borrowed $1$ and divides $9.0.4.10.4$ by $2$
giving him $a = 10.0.2.5.2.7$ instead of the correct value.

The $11$-adic numbers $a$ and $b$ given by Gauss satisfy $a + b = 10$ and
$a - b = \sqrt{5}$. Since there is an error in the values of $a$
and $b$ it is difficult to say what Gauss is doing here; at any rate,
the fourth number is the sum of the third and the fifth.

\section{Square roots of $1$ in $10$-adic numbers}

Today, $p$-adic numbers are everywhere. In contrast, $g$-adic
numbers for composite integers $g$ do not play a major role in number
theory. From an algebraic point of view, $\Z_{pq} \simeq \Z_p \oplus \Z_q$
is just the direct sum of two $p$-adic rings whenever $p$ and $q$ are
distinct primes; in particular, $\Z_{pq}$ has zero divisors.

Let us now see what's behind the isomorphism $\Z_{10} \simeq \Z_2 \oplus \Z_5$.
There are natural projections
$\Z_{10} \lra \Z_2$ and $\Z_{10} \lra \Z_5$ induced by sending a residue class
modulo $10^n$ to the residue classes modulo $2^n$ and $5^n$. The Chinese
Remainder Theorem  yields
isomorphisms $\Z/10^n\Z \simeq \Z/2^n\Z \oplus \Z/5^n\Z$, which then
implies the desired isomorphism by taking limits.

$10$-adic numbers do, however, show up in recreational introductions to
``unusual'' number systems\footnote{See e.g. Borcherdt's introduction
  \cite{Borch}.}. A favourite topic in this area is the
convergence of numbers whose square has the same end digits,
namely\footnote{Such problems were discussed
  starting in 1815; see \cite{Gerg} and \cite{Schreven}.
  Paul Z\"uhlke (a teacher of mathematics in Prussia) \cite{Zuehlke}
  gives several references to publications in the 19th century and remarks
  that Hensel has communicated a short (but not elementary) solution of
  the congruence $x^2 \equiv x \bmod 10^n$. Koppe \cite{Koppe} employed
  continued fractions.}
$$ a = \ldots918\,212\,890\,625. $$
In fact, $890\,625^2 = 793\,212\,890.625$ and
$2\,890\,625 = 8\,355\,712\,890\,625$ etc. 
This $10$-adic number satisfies the equation $a^2 = a$, hence
$0 = a^2 - a = a(a-1)$, which implies that $\Z_{10}$ has zero divisors.

The numbers $a$ and $b = 1-a$ 
\begin{align*}
    a & = \ldots 392\,256\,259\,918\,212\,890\,625 \quad \text{and}  \\
    b & = \ldots 607\,743\,740\,081\,787\,109\,376
\end{align*}
are both ``automorphic'', i.e., we have $a^2 = a$, which in turn implies
that we also have $b^2 = (1-a)^2 = 1 - 2a + a^2 = 1 - 2a + a = 1-a = b$.

\subsection*{The algebra behind these numbers}

Let us now briefly explain the algebraic background of the two
automorphic numbers $a$ and $b$. This is, as far as I can see, not
directly relevant for Gauss's calculations, but it helps us understand 
how {\tt pari} verifies the results.

The elements $a$ and $b = 1-a$ are idempotents in the ring of
$10$-adic integers, which means that $a^2 = a$ and $b^2 = b$. In addition,
they are orthogonal, i.e., $ab = 0$. Moreover, the unit $1 = a + b$ is a
sum of orthogonal idempotents. Given any ring in which $1$ is the sum of
two orthogonal idempotents $a$ and $b$, we can write every element $r \in R$
as $r = r \cdot 1 = r (a + b) = ra + rb$. It is easy to see that this
allows us to write the ring $R$ as a direct sum $R = Ra \oplus Rb$.

The pre-images $a$ and $b$ of the elements $(1,0)$ and $(0,1)$ in
$\Z_2 \oplus \Z_5$ with respect to the projection $\Z_{10} \lra \Z_2 \oplus \Z_5$
are idempotents since $(1,0) \cdot (1,0) = (1,0)$, 
and  they are orthogonal since $(1,0) \cdot (0,1) = (0,0)$.
The elements $a$ and $b$ can easily be computed using {\tt pari} via
\begin{align*}
  a & = {\tt chinese(Mod(0,5^{20}),Mod(1,2^{20}))}, \\
  b & = {\tt chinese(Mod(1,5^{20}),Mod(0,2^{20}))}.
\end{align*}

We can also show that $a = \lim 5^{2^n}$; the sequence $a_n = 5^{2^n}$
converges $10$-adically since $a_{n+1} - a_n \equiv 0 \bmod 10$. Similarly,
we have $b = \lim 6^{5^n}$.

$$ \begin{array}{r|r|r}
  \rsp n & 5^{2^n} & 6^{5^n} \\ \hline
  \rsp 1 &           25 &  7776 \\
  \rsp 2 &          625 &  \ldots 1376 \\
  \rsp 3 & \ldots 90625 &  \ldots 9376
  \end{array} $$

The two decompositions of $\Z_{10}$ as direct sums, namely
$\Z_{10} \simeq \Z_{2} \oplus \Z_{5}$ and
$\Z_{10} \simeq \Z_{10}a \oplus \Z_{10}b$, are actually the same: it
is clear that the homomorphism $\pi_2: \Z_{10} \lra \Z_{2}$
defined by reduction modulo large powers of $2$ satisfies
$\pi_2(\Z_{10}a) = \Z_2$ since $\pi_2(a) = 1$ and $\pi_5(a) = 0$.

\subsection*{Square roots of $1$}

If $i$ is an idempotent, then
$$ (2i-1)^2 = 4i^2 - 4i + 1 =  4i - 4i + 1 = 1, $$
hence $2a - 1$ is a square root of $1$ in $R$ different from $\pm 1$.
Similarly, $2b - 1 = 2(1-a)-1 = -(2a-1)$ is a square root of $1$.
These two square roots of $1$ given  are the pre-images of $(-1,+1)$
and $(+1,-1)$ in $\Z_2 \oplus \Z_5$, and pari shows that
\begin{align*}
    2a-1 & = {\tt chinese(Mod(-1,5^{20}),Mod(+1,2^{20}))}
           = \ldots 8\,451\,251\,9836\,425\,781\,249, \\
    2b-1 & = {\tt chinese(Mod(+1,5^{20}),Mod(-1,2^{20}))}
           = \ldots 1\,5487\,480\,163\,574\,218\,751.
\end{align*}

On p.~40 of \cite{GS}, Gauss computes the square root $\eps = 2a-1$ of $1$
modulo $10^\infty$; he finds
$$ \eps = \ldots 425\,781\,249. $$
The calculation of the square root using the binomial expansion is not
very user-friendly since $\eps$ is not an integer, which is why Gauss
adapts the algorithm for computing square roots digit by digit.

\subsection*{Square Root Algorithm in the $p$-adic Numbers}

Gauss uses the algorithm for computing square roots digit by digit
(from right to left) for computing the $10$-adic square root $\eps$
of $1$ with $\eps \equiv 49 \bmod 100$. In the $10$-adics we cannot
use the size of the numbers; the new digits must be determined using
congruences.

We start by subtracting $49^2 = 2401$ from $1$:

\begin{center}
  \begin{tikzpicture}
    \draw (0,0.5) node[anchor=east] {$49$}; 
    \draw (0,0) node[anchor=east] {$\sqrt{\ldots |00|00|01}$};
    \draw (-3,0) node[anchor=east] {$49$};
    \draw (0,-0.4) node[anchor=east] {$24|01$};
    \draw (0,-0.8) node[anchor=east] {$\ldots99|76|00$};
  \end{tikzpicture}
\end{center}

\subsubsection*{The algorithm}
We will now present Gauss's algorithm by explaining how to proceed
once the new digit has been found, then we show how to find the digit,
and finally we explain why the algorithm works.

\begin{enumerate}
\item If the remainder is $\equiv r \bmod 100$ for some $2$-digit
  number $r$, then the new digit is $ \equiv - \frac r2 \bmod 10$.
\item Double the number at the top and append the new digit $n$ on the left.
  The result is written down on the left, the new digit is written on top.
\item Multiply the number on the left by the new digit and subtract.
\end{enumerate}
With $r = 76$, the first new digit is $b \equiv - 38 \equiv 2 \bmod 10$,
hence $b = 2$. The next new digits are $1$, $2$, and $8$. The calculation
of the square root thus looks as follows:

\begin{center}
  \begin{tikzpicture}
    \draw (0,0.5) node[anchor=east] {$81249$}; 
    \draw (0,0) node[anchor=east] {$\sqrt{\ldots |00|00|01}$};
    \draw (-3,0) node[anchor=east] {$49$};
    \draw (0,-0.4) node[anchor=east] {$24|01$};
    \draw (0,-0.8) node[anchor=east] {$\ldots 99|76\phantom{|00}$};
    \draw (0,-1.2) node[anchor=east] {$\ldots  5|96\phantom{|00}$};
    \draw (-3,-0.8) node[anchor=east] {$298$};
    \draw (0,-1.6) node[anchor=east] {$\ldots 9|93|8\phantom{0|00}$};
    \draw (-3,-1.6) node[anchor=east] {$1498$};
    \draw (0,-2) node[anchor=east] {$ 1|49|8\phantom{0|00}$};
    \draw (-3,-2.4) node[anchor=east] {$82498$};
    \draw (0,-2.4) node[anchor=east] {$\ldots 99|98|44\phantom{|00|00}$};
    \draw (0,-2.8) node[anchor=east] {$ 65|99|84\phantom{|00|00}$};
    \draw (0,-3.2) node[anchor=east] {$\ldots 33|98|6\phantom{0|00|00}$};    
  \end{tikzpicture}
\end{center}

Doubling $81\,249$ and omitting the digit $1$ on the left, we find
the new digit $7$, hence we have to form
$762\,498 \cdot 7 \equiv 37486 \bmod 10^6$. Gauss continues and finds the
next three digits  $5$, $2$ and $4$, giving him
$$\eps = \ldots 425\,781\,249. $$

\subsubsection*{Why it works}
Gauss starts with the $10$-adic approximation $\eps \equiv 49 \bmod 100$.
For finding the next digit $a$, we compute
$$ 1 - (100a+49)^2 = \ldots 9\,997\,600 - 9\,800a - 10^4a^2, $$
or, if we discard the last two zeros as Gauss does, then
$$ 0 \equiv  99\,976 - 98a - 100a^2  \equiv 99\,976 - 98a \bmod 100. $$
This implies $49a \equiv 88 \bmod 50$ and thus $a \equiv 2 \bmod 10$.
In particular, the last three digits of $\eps$ are $249$.


\begin{figure}[ht!]
  \includegraphics[width=4cm]{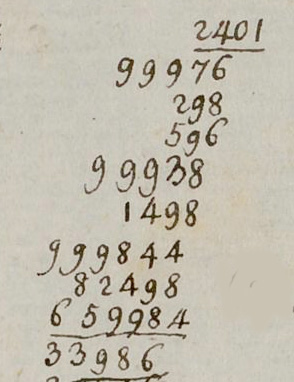}
  \caption{Calculation of a $10$-adic square root of $1$}
\end{figure}

Now Gauss has to subtract 
$$ 400  + 2 \cdot 98 = 2 \cdot 298 , $$
and finds that $10^\infty + 1 - 249^2  = \ldots 99\,938$. From
$$ 10^\infty + 1 - (1000b + 249)^2
   =  \ldots 99\,938\,000 - 2 \cdot 249\,000b - 10^6 b^2 $$
we now get
$$  98b \equiv 38 \bmod 100, $$
which implies $b = 1$ and
$$  10^\infty + 1 - (1249)^2 = 99\,938\,000 - 498\,000-1\,000\,000, $$
or $99\,938 - 1\,498 = 98\,440$.

In the general case, assume that $1 - a^2 \equiv 10^nr \bmod 10^{n+2}$,
so $r$ is a two-digit integer (whose last digit is even by construction).
Setting
$$ 0 \equiv 1 - (a+10^nb)^2 \bmod 10^{n+2} $$
yields
$$ 10^nr  \equiv - 2ab\cdot 10^n \bmod 10^{n+2}, $$
that is,
$$ b \equiv - \frac{r}{2} \bmod 10 $$
since $a \equiv 9 \bmod 10$.

From $1 - 249^2 = \ldots 99\,938\,000$ we get $r = 38$ and
$b \equiv - 19 \equiv 1 \bmod 10$. 
Now $1 - 1249^2 = \ldots 999\,844$, hence $r = 44$ and
$b \equiv -22 \equiv 8 \bmod 10$.
Continuing in this way we find more digits of this square root of $1$.

Observe that if $n$ is large enough, one may discard the term
$10^{2n} \cdot b^2$ and just work with the last nonzero digits, which is
what Gauss is doing in the last few lines of his calculation.

\section{Logarithms in $10$-adics}

Gauss's calculations of the $10$-adic square roots of $1$ is on the
same page as his calculations of $10$-adic logarithms of natural numbers.
This might be a coincidence -- I am not aware of any connection between
the $10$-adic square roots of $1$ and logarithms.

\subsection*{Calculation of $\log 2$}
Gauss starts by computing $\log(2) = \ldots 21.830.960$; I do not know
how he did it. Since his calculation is just below the value of $\eps$
I suspect that there may be a connection.

For us it seems natural to use the isomorphism
$\Z_{10} \simeq \Z_2 \oplus \Z_5$ and then define $\log(2)$ as
the preimage of $(0,\log_5(2))$. This yields $\log(2) = \ldots863\,080\,960$,
which agrees with Gauss's $\log(2)$ only modulo $10^4$.
Observe that
$$ \log_5(2) \equiv 34\,085 \bmod 5^7 \quad \text{and} \quad
   21\,830\,960 \equiv 34\,085 \bmod 5^7. $$
On the other hand, 
$$ 21\,830\,960 \equiv 48 \bmod 2^7. $$
Thus it seems that Gauss made an error in his calculations.

For computing $\log_5(2)$ we can use (\ref{logx-1}) with $x = 5$ and find
$$ -2 \log_5(2) = -\log_5(4) = 5 + \frac{5^2}2 + \frac{5^3}3 + \ldots, $$
which gives $\log_5(2) \equiv 335 \bmod 5^4$ and now
$$ 335 = 2 \cdot 5 + 3 \cdot 5^2 + 2 \cdot 5^3 $$
gives us the first three digits of 
$$ \log_5(2) = 
   2 \cdot 5 + 3 \cdot 5^2 + 2 \cdot 5^3 + 4 \cdot 5^4 +
   2 \cdot 5^6 + 2 \cdot 5^7 + 4 \cdot 5^8 + \ldots. $$
The $5$-adic approximation of $\log(2)$, namely
$$ 2 \cdot 5 + 3 \cdot 5^2 + 2 \cdot 5^3 + 4 \cdot 5^4 +
   2 \cdot 5^6 + 2 \cdot 5^7 + 4 \cdot 5^8 + 2 \cdot 5^9 = 5659085 $$
can be transformed into a $10$-adic number by multiplication
with the idempotent $b = \ldots787\,109\,376$, and in fact
$$ 5\,659\,085 \cdot 787\,109\,376 \equiv 886\,308\,0960 \bmod 10^{10}. $$
I do not think that Gauss proceeded this way.   

\subsection*{Calculation of $10$-adic logarithms using power series}

The power series
\begin{equation}\label{logx+1}
    \log(1+x) = x - \frac{x^2}2 + \frac{x^3}3 - \frac{x^4}4 + \ldots 
\end{equation}
converges for integers $x \equiv 0 \bmod 10$; for example,
$$ \log (31) \equiv 30 - \frac{30^2}2 + \frac{30^3}3 - \frac{30^4}4
    + \frac{30^5}5-\frac{30^6}6 - \frac{30^8}{8} \ldots
    \equiv 666\,080 \bmod 10^7 $$
since $\frac17 \cdot 30^7 \equiv 0 \bmod 10^7$,     
and the value Gauss gives is $\log(31) = \ldots 80\,666\,080$.
The logarithm of numbers coprime to $10$ can be computed using the
power series; for example, $4\log(3) = \log(81)$, and the result
agrees with Gauss's calculation.

\begin{figure}[ht!]
  \begin{center}
    \includegraphics[height=4cm]{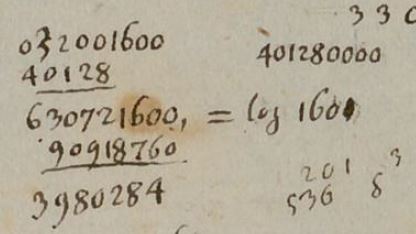}
  \end{center}
  \caption{Gauss's calculation of $\log(1601)$.}\label{1601}
\end{figure}

Gauss did use the power series, as the calculations in the left
picture of Fig.~\ref{1601} show.
\begin{align*}
  \log (1601) & \equiv 1600 - \frac12 1600^2 + \frac13 \cdot 1600^3- \ldots \\
       & \equiv 1600 - 1\,280\,000 + 32\,000\,000 - 400\,000\,000 + \ldots \\
       & \equiv  32\,001\,600 - 401\,280\,000 \equiv 630\,721\,600 \bmod 10^9.
\end{align*}
{\tt pari} gives $\log(1601) = \ldots 1\,630\,721\,600$.

At the bottom of the page he calculates $\log(59)$
using the power series
\begin{equation}\label{logxm}
  - \log(x-1) = x + \frac{x^2}2 + \frac{x^3}3 + \frac{x^4}4 + \ldots.
\end{equation}With $x = 60$ he finds
\begin{align*}
  - \log 59
  & = 60 + \frac{60^2}2 + \frac{60^3}3 + \frac{60^4}4 + \frac{60^5}5
      + \frac{60^6}6 + \ldots \\ 
  & = 60 + 1\,800 + 72\,000 + 3\,240\,000 + 155\, 520\,000
      + 7\,776\,000\,000 + \ldots \\
  & = \ldots 34833860
\end{align*}
hence
$$ \log(59) = - \ldots 34833860 = \ldots 85\,166\,140. $$

\begin{figure}[ht!]
  \begin{center}
    \includegraphics[height=4cm]{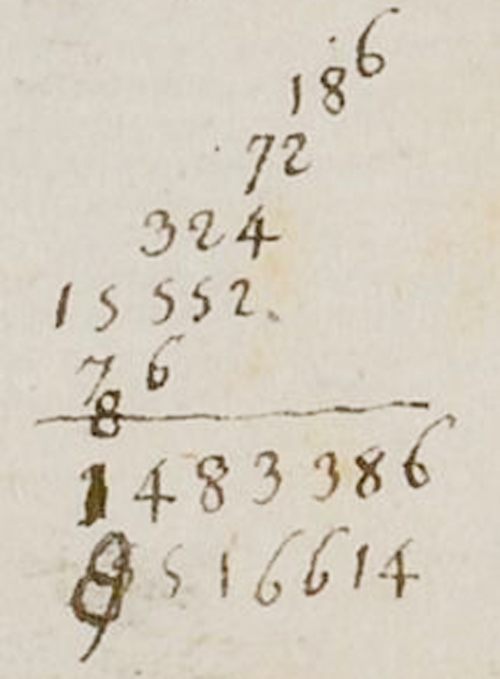} \quad 
    \includegraphics[height=4cm]{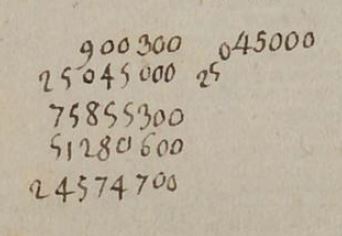} \quad 
  \end{center}
  \caption{Gauss's calculation of $\log(59)$ (left) and
    $\log(301)$ (right).}\label{log59}
\end{figure}

His calculation of $\log (301)$ contains an error: we have
\begin{align*}
  \log(301) & \equiv 300 - \frac12 \cdot 300^2 + \frac13 \cdot 300^3
               - \frac14 \cdot 300^4 + \ldots \\
            & \equiv 300 - 45\,000 + 9\,000\,000 - 25\,000\,000 + \ldots \\
            & \equiv 9\,000\,300 - 25\,045\,000 \equiv 983\,955\,300 \bmod 10^9, 
\end{align*}
but Gauss omits a $0$ in $\frac12 \cdot 300^2 = 9\,000\,000$ and gets
$75\,855\,300 \bmod 10^8$ (Fig.~\ref{1601}, right).

The numbers at the bottom are
$$ \log( 7) = \ldots 51\,280\,600 \quad \text{and} \quad
   \log(43) = \ldots 32\,674\,700 $$
whose sum is $\log(301)$ since $301 = 7 \cdot 43$.
Gauss's table of $10$-adic logarithms is complete up to $25$;
then he only gives the logarithms of a few primes, namely
$29$, $31$, $43$ and $59$.

\begin{center}
  \includegraphics[height=4cm]{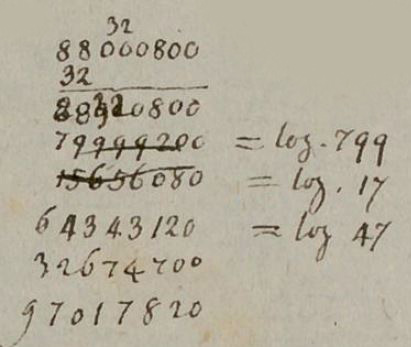} 
\end{center}

Here Gauss computed $\log(799)$ by setting $x = 800$ in the power series
\begin{equation}\label{logx-1}
  - \log(x-1) = x + \frac{x^2}2 + \frac{x^3}3 + \frac{x^4}4 + \ldots
\end{equation}
and he finds
$$ -\log 799 \equiv 800 + 320\,000 \equiv 320\,800 \bmod 10^5. $$
The next term would be $\frac13 \cdot 800^3 \equiv 4\,000\,000 \bmod 10^8$;
again there seems to be an error in Gauss's computations. The correct
values are
\begin{align*}
  \log(799) & = \ldots 95\,679\,200, \\
  \log(17)  & = \ldots 15\,656\,080, \\
  \log(47)  & = \ldots 80\,023\,120.
\end{align*}

\begin{center}
  \includegraphics[height=4cm]{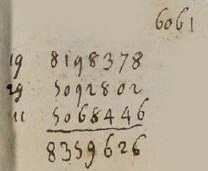} 
\end{center}

Another calculation by Gauss is the logarithm of
$6061 = 11 \cdot 19 \cdot 29$. The correct values are
$$ \log(19) = \ldots 81\,983\,780,
   \log(29) = \ldots 50\,928\,020 \quad \text{and} \quad 
   \log(11) = \ldots 50\,684\,460 $$
giving $\log(6061) = \ldots 83\,596\,260$.
   
{\tt pari} can compute $10$-adic logarithms using the isomorphism
$\Z_{10} \simeq \Z_2 \oplus \Z_5$. We compute the $2$-adic and $5$-adic
logarithm and then the Chinese Remainder Theorem will give us the
$10$-adic logarithm. In our example,
$$ {\tt n=12;
    chinese(Mod(lift(log(31+O(2^n))),2^n),Mod(lift(log(31+O(5^n))),5^n))} $$
yields
$$ \log(31) = \ldots 723\,280\,666\,080. $$

\begin{figure}[ht!]
  \includegraphics[height=4cm]{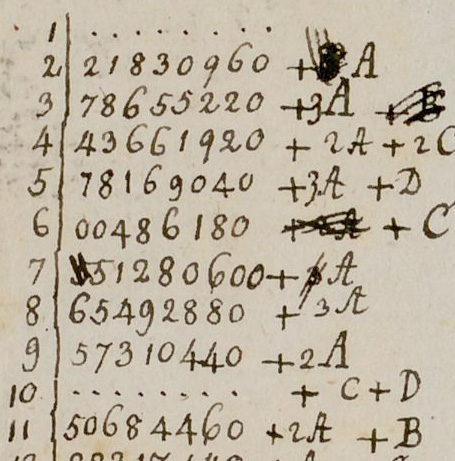} \quad 
  \includegraphics[height=4cm]{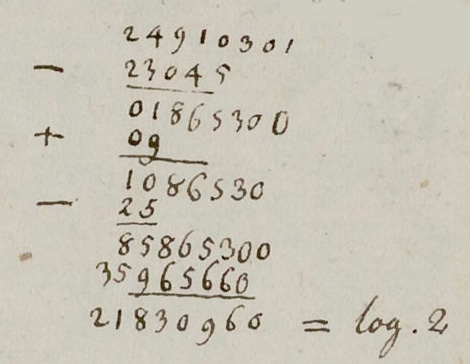}
  \caption{Calculation of $10$-adic approximations of logarithms}\label{Fl2}
\end{figure}

The logarithms of the integers Gauss has computed satisfy
$\log(ab) = \log(a) + \log(b)$. He does not give $\log(1)$ and $\log(10)$,
but the latter can be computed by adding $\log(2)$ and $\log(5)$; the result
is $\log(10) = 0$. This is in accordance with Gauss's result
$\log(20) = \log(2)$.

On the second page of his calculations Gauss observes that
$$ \log(a+b) = \log(a) + \frac ba - \frac{b^2}{2a^2} + \ldots. $$
In fact,
\begin{align*}
  \log(a+b) - \log(a)
  & = \log\Big(\frac{a+b}a\Big)
    = \log \Big(1 + \frac ba\Big) \\
    & = \frac ba - \frac{b^2}{2a^2} + \frac{b^3}{3a^3} - \frac{b^4}{4a^4} + \ldots
\end{align*}
Here we need $|\frac ba| < 1$ in the real numbers, and 
$10a \mid b$ for convergence in the $10$-adic numbers.

\begin{figure}[ht!]
  $$ \begin{array}{r|r}    
    1 & 0 \\
    2 & 63080960 \\
    3 & 78655220 \\
    4 & 26161920 \\
    5 & 98437500 \\
    6 & 91736180 \\
    7 & 51280600 \\
    8 & 89242880 \\
    9 & 57310440 \\
    10 & 21518460
  \end{array} \quad 
  \begin{array}{r|r}
    11 & 50684460 \\
    12 & 29817140 \\
    13 & 97535940 \\
    14 & 64361560 \\
    15 & 57092720 \\    
    16 & 52323840 \\
    17 & 15656080 \\
    18 & 70391400 \\    
    19 & 81983780 \\
    20 & 19599420
  \end{array} \quad 
  \begin{array}{r|r}
    21 & 29935820  \\    
    23 & 68418760 \\
    25 & 96875000 \\
    27 & 35965660 \\
    29 & 50928020 \\
    31 & 80666080 \\
    33 & 29339680 \\
    35 & 29718100 \\
    37 & 63533340 \\
    39 & 76191160
  \end{array} $$  
  \caption{Correct table of 10-adic logarithms}
\end{figure}

\begin{figure}[ht!]
  \includegraphics[height=4cm]{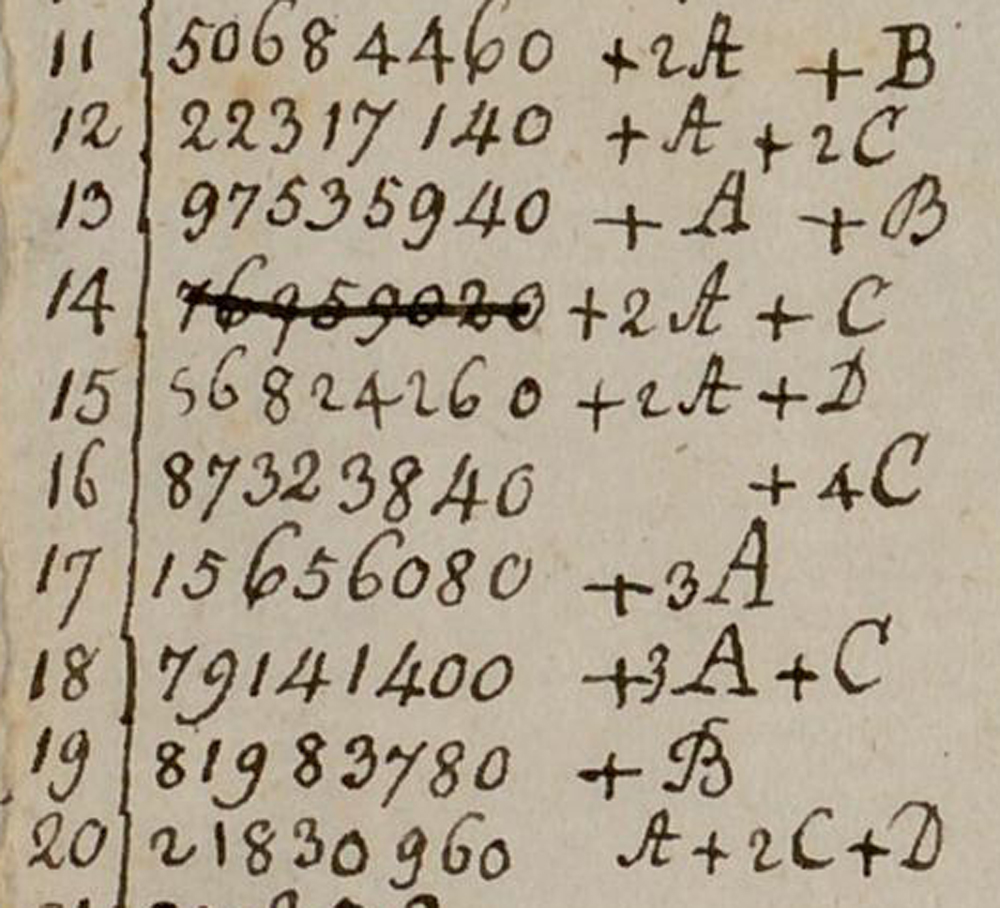} \quad
  \includegraphics[height=4cm]{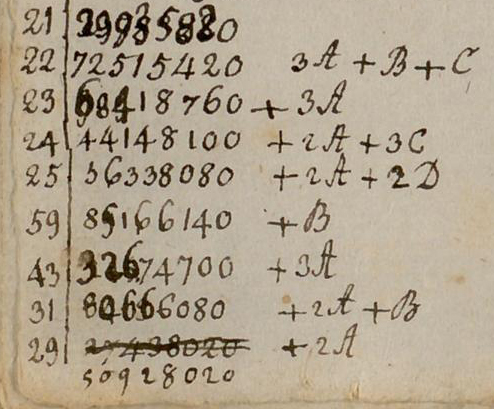}  
  \caption{More logarithms from Gauss's notebooks}
\end{figure}

\begin{figure}[ht!]
  \includegraphics[width=3.29cm]{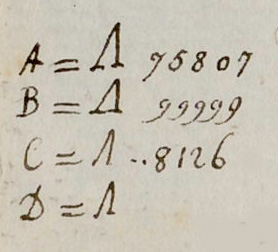}
  \caption{Four $10$-adic numbers}
\end{figure}

Gauss also defines four numbers
\begin{align*}
  A & = \Lambda 75807, \\
  B & = \Lambda 99999, \\
  C & = \Lambda..8126, \\
  D & = \Lambda
\end{align*}
whose meaning is not clear to me. As elements in $\Z_2 \oplus \Z_5$ it
seems that
$$ \begin{array}{clll}
  A & = (-1,\sqrt{-1}) & = \ldots666\,295\,807, & 
      \sqrt{-1} = 2 + 5 + 2 \cdot 5^2 + 5^3 + 3 \cdot 5^4 + \ldots  \\
  B & = (-1,-1) & = \ldots999\,999\,999, &  \\
  C & = (-2, 1) & = \ldots361\,328\,126. & 
  \end{array} $$
Observe that $A^2 = (1,-1)$ = , so $A$ is a square root of $\eps$.
These numbers occur next to the logarithms of the natural numbers in
Fig.~\ref{Fl2}, but I do not understand their meaning.

I observe that the logarithms of numbers $\equiv \pm 3 \bmod 20$
have the entry $+ 3A$, and numbers $\equiv 9 \bmod 20$ have $+ 2A$.

\section*{Comments}
In his article {\em The unpublished Section Eight: On the way to
  function fields over a finite field} (see \cite[Ch. IV]{GSS}),
G\"unther Frei pointed out that
Gauss was in possession of ``Hensel's Lemma'' for lifting roots of a
polynomial modulo a prime number $p$ to roots modulo higher powers of $p$.
His notebook reveals that, at the same time, he knew how to do basic
arithmetic with ``infinite congruences''; these infinite congruences
coincide with the modern idea of compatible systems of congruences
(inverse limits). In his notes, Gauss did not formally define these
limits, but he knew how to add, subtract and multiply them. Gauss also
could compute square roots, define a $p$-adic logarithm using the Taylor
expansion of $\log(1+x)$, and extend it beyond its domain of convergence.
More detective work is required to decrypt all of Gauss's calculations
concerning $10$-adic logarithms.

After having written the first version of this note I have discovered
that Gauss's infinite congruences have been noticed by user2554
(see also \cite{Conrad}) in a posting on math stackexchange \cite{MSE}
dated March 8, 2024.

\end{document}